\documentstyle[12pt]{article}
\newcommand{\const}{\mathop{\rm const}\limits}

\newcommand{\grad}{\mathop{\rm grad}\limits}

\newcommand{\supp}{\mathop{\rm supp}\limits}

\begin{document}
\begin{center}

{\bf Module of continuity for the functions }\\
\vspace{2mm}
{\bf  belonging  to the Sobolev-Grand Lebesgue Spaces }\\

\vspace{3mm}

{\sc Ostrovsky E., Sirota L.}\\


\vspace{3mm}
 Department of Mathematics and Statistics, Bar-Ilan University,
59200, Ramat Gan, Israel.\\
e-mails: \ galo@list.ru, \  sirota@zahav.net.il \\

\end{center}

\vspace{3mm}

 {\it Abstract.} In this short article we generalize the Sobolev's  inequalities
 for the module of continuity for the functions belonging to the classical Lebesgue
space on the (Bilateral) Grand Lebesgue spaces. \par
 We construct also some examples in order to show the exactness of obtained results.\par
 \vspace{3mm}

  \par

 \vspace{3mm}


\vspace{3mm}

\section{Introduction. Notations. Statement of problem.}

\vspace{3mm}
{\bf 1. Sobolev's continuity inequality.} \par
Let $ D $ be convex closed bounded non-empty domain with Lipshitz boundary
in the Euclidean space $ R^d, \ d \ge 2. $ For instance, the domain $ D $ may be the
unit ball $ B $ on the space $ R^d: $

$$
B = \{ x: x \in R^d, \ |x| \le 1 \}.
$$

We will consider in this article only the case when the domain $ D $ is bounded.\par
 The classical Sobolev's inequality (for the domain $ D $ or for
 the whole space $ R^d,) $ see, e.g.
 \cite{Kantorovicz1}, chapter 11, section 5; \cite{Sobolev1}, \cite{Talenti1} etc.
 asserts that for all weak differentiable
 functions $ f, f: R^d \to R, \ d \ge 3 $ from the Sobolev's space
 $ W^1_p(D), \ p \in [1,d), $ which may be defined as a closure in the Sobolev's norm

 $$
 ||f||W^1_p(R^m)= |f|_p + |\nabla f|_p
 $$
  of the set of all finite continuous differentiable functions $ f, f: D \to R, $ that

  $$
  |f|_{q} \le K_d(p) \ |\nabla f|_p, \ q = q(p)= dp/(d-p), \ p \in [1,d), \
  q \in (d/(d-1), \infty). \eqno(1)
  $$
    Here the notation $ |x| $ denotes ordinary Euclidean norm of the vector $ x, $

  $$
  |f|_p = |f|_{p,D} = \left[ \int_{D} |f(x)|^p \ dx  \right]^{1/p},
  $$

  $$
  \nabla f = \{ \partial f/\partial x_1, \partial f/\partial x_2, \partial f/\partial x_3, \ldots,  \partial f/\partial x_d \} = \grad_x f,
  $$

$$
|\nabla f|_p = \left| \left[ \sum_{i=1}^d ( \partial f/\partial x_i)^2  )  \right]^{1/2} \right|_p,
$$

$$
|\nabla^2 f|_p =
\left| \left[ \sum_{i=1}^d \sum_{j=1}^d ( \partial^2 f/\partial x_i \partial x_j)^2)  \right]^{1/2} \right|_p
$$
etc.\par
The best possible constant in the inequality (1) belongs to G.Talenti \cite{Talenti1}:

$$
K_d(p) = \pi^{-1/2} d^{-1/p} \left[\frac{p-1}{d-p}\right]^{1-1/p} \cdot
\left[\frac{\Gamma(1+d/2) \ \Gamma(d)}{\Gamma(d/p) \ \Gamma(1+d-d/p)} \right]^{1/d}.
$$
The case $ p = d $ was considered by Yudovich \cite{Yudovich1} and after 7 years by
Trudinger \cite{Trudinger1}; see also \cite{Saloff-Coste1},section 1,2;
\cite{Taylor1}, chapter 13, section 4.\par
The exact values of constants  in the Orlich-Sobolev imbedding theorem in this case
$ p = d $ was obtained by Moser \cite{Moser1}.\par
\vspace{3mm}
 Let us consider hereafter {\it only} the case $ p > d, \ d = 1,2,\ldots; $
 and also (for simplicity) {\it only } the case when
 $$
f(x)_{\partial D} \stackrel{def}{=} \lim_{x \to \partial D} f(x) = 0;
 $$
 where $ \partial D $ denotes the boundary of the domain $ D.$ \par
 In this case, i.e. when
$ \nabla f \in L_p(D), $ the function $ f(\cdot) $ up to chaining in the subset of
the domain $ D $ of the zero measure is continuous. \par
 More detail, let us define as ordinary the module of continuity $ \omega(f,\delta) $
 of arbitrary uniform continuous function $ f: D \to R $

$$
\omega(f,\delta) \stackrel{def}{=} \sup_{x,y \in D, |x-y| \le \delta} |f(x) - f(y)|,
$$
$ \delta \in (0,1/e).$ \par
It is proved, e.g. in the book \cite{Maz'ja1}, chapter 1, p. 60-62  that
$$
\omega(f,\delta) \le K_M(d,D) \ \delta^{1-d/p} \cdot
\left[ \frac{p-1}{p-d} \right]^{1-1/p} \cdot |\nabla f|_p. \eqno(2a)
$$
Analogous result is obtained by Leoni \cite{Leoni1}, chapter 11, section 11.3 with the
exact values of a constant:

$$
\omega(f,\delta) \le \frac{2dp}{p-d} \cdot (2 \delta)^{1-d/p} \cdot |\nabla f|_p.\eqno(2b)
$$
Some important applications of these inequalities in the theory of Partial Differential
Equations are described, e.g., in \cite{Evans1}, \cite{Taylor1}; in the Calculus of
Variations - in \cite{Morrey1}. \par
 Notice that both the inequalities may be rewritten in the equivalent up to
 multiplicative constant form:

$$
\omega(f,\delta) \le C(d,D) \frac{p}{p-d} \cdot \delta^{1-d/p} \cdot |\nabla f|_p.\eqno(2c),
$$
as long as $ p > d. $ \par
In the book  \cite{Evans1}, chapter 5, section 5.6, p. 280-282
the inequalities (2a), (2b) and (2c) was named as a
Morrey inequality; see also \cite{Morrey1}. \par

 The inequalities (2a), (2b) and (2c) was generalized in the works \cite{Adams1},
 \cite{Cianchi1}, \cite{Donaldson1} on the Orlicz-Sobolev spaces, i.e.
 when $ \nabla f $ belongs to some Orlicz space. We intent to improve, in particular,
 this results. \par

\vspace{3mm}

{\bf 2. Our aim.}\par
{\bf Our aim is generalization of Sobolev's continuity inequality (2a), (2b) or (2c)
on some popular classes of rearrangement invariant (r.i.) spaces, namely, on the so-called Sobolev-Grand Lebesgue Spaces $ G(\psi). $  We intend to show also the exactness of offered estimations.}

 \vspace{3mm}

{\bf 3. Another notations.}
 Hereafter $ C, C_j $ will denote any non-essential finite positive constants.
 We define also for the values $ (p_1, p_2), $ where $ 1 \le p_1 < p_2 \le \infty $

$$
L(p_1, p_2) =  \cap_{p \in (p_1, p_2)} \ L_p.
$$

 We denote also

 $$
 \Omega(d)= \frac{2 \pi^{d/2}}{\Gamma(d/2)};
 $$
and denote as usually an indicator function

  $$
 I(A)= I(A,x) = 1, x \in A, \ I(A) = I(A,x) = 0, \ x \notin A.
 $$

\vspace{3mm}

{\bf 4. Content of the paper. } \par
\vspace{3mm}

  The paper is organized as follows. In the next section we recall the definition
and some simple properties of the so-called Grand Lebesgue Spaces (GLS)
$ G(\psi) $ and introduce its generalizations: the so-called Sobolev's Grand Lebesgue
spaces (SGLS) $ W^{1}G(\psi). $ \par
 In the section 3 we formulate and prove the main result: the classical Sobolev's
 continuity inequality  for $ W^{1}G(\psi) $ spaces. \par
 In the fourth section we built some examples in order to show the exactness of obtained
 inequalities. The section 5 is devoted of the consideration the one-dimensional case,
 in which we can compute the exact value on embedding constant.\par

  The last section contains some concluding remarks: a {\it weight} generalizations of
  the embedding theorem for the hight derivatives etc. \par

\vspace{4mm}

\section{Sobolev's-Grand Lebesgue Spaces.}\par

 {\bf A. Grand Lebesgue Spaces.} \par

\vspace{3mm}

     Recently, see \cite{Kozachenko1}, \cite{Fiorenza1},
     \cite{Fiorenza2}, \cite{Fiorenza3}, \cite{Iwaniec1}, \cite{Iwaniec2},
     \cite{Ostrovsky1}, \cite{Ostrovsky2}, \cite{Ostrovsky3}, \cite{Ostrovsky4},
     \cite{Ostrovsky5}, \cite{Ostrovsky6}, \cite{Ostrovsky8}, \cite{Ostrovsky9}  etc.
     appears the so-called Grand Lebesgue Spaces $ GLS = G(\psi) =G\psi =
    G(\psi; A,B), \ A,B = \const, A \ge 1, A < B \le \infty, $ spaces consisting
    on all the measurable functions $ f: T \to R $ with finite norms

     $$
     ||f||G(\psi) \stackrel{def}{=} \sup_{p \in (A,B)} \left[ |f|_p /\psi(p) \right].
     \eqno(3)
     $$

      Here $ \psi(\cdot) $ is some continuous positive on the {\it open} interval
    $ (A,B) $ function such that

     $$
     \inf_{p \in (A,B)} \psi(p) > 0, \ \psi(p) = \infty, \ p \notin (A,B).
     $$
We will denote
$$
 \supp (\psi) \stackrel{def}{=} (A,B) = \{p: \psi(p) < \infty, \} \eqno(4)
$$

The set of all $ \psi $  functions with support $ \supp (\psi)= (A,B) $ will be
denoted by $ \Psi(A,B). $ \par
  This spaces are rearrangement invariant, see \cite{Bennet1}, and
  are used, for example, in the theory of probability \cite{Talagrand1}, \cite{Kozachenko1}, \cite{Ostrovsky1}; theory of Partial Differential Equations \cite{Fiorenza2}, \cite{Iwaniec2}; functional analysis \cite{Ostrovsky4}, \cite{Ostrovsky5}; theory of Fourier series \cite{Ostrovsky7}, theory of martingales \cite{Ostrovsky2} etc.\par
 Notice that in the case when $ \psi(\cdot) \in \Psi(A,B),  $ a function
 $ p \to p \cdot \log \psi(p) $ is convex, and  $ B = \infty, $ then the space
$ G\psi $ coincides with some {\it exponential} Orlicz space. \par
 Conversely, if $ B < \infty, $ then the space $ G\psi(A,B) $ does  not coincides with
 the classical rearrangement invariant spaces: Orlicz, Lorentz, Marzinkievitch etc.
\cite{Ostrovsky8}, \cite{Ostrovsky9}. \par
 We will use the following two important examples (more exact, the {\it two families
of examples} of the $ \psi $ functions and correspondingly the GLS spaces.\par
{\bf 1.} We denote
$$
\psi(A,B; \alpha,\beta;p) \stackrel{def}{=} (p-A)^{-\alpha} \ (B - p)^{-\beta},\eqno(5)
$$
where $ \alpha,\beta = \const \ge 0, 1 \le A < B < \infty; p \in (A,B) $ so that
$$
\supp \psi(A,B; \alpha,\beta;\cdot) = (A,B).
$$

{\bf 2.} Second example:
$$
\psi(1,\infty; 0, -\beta;p) \stackrel{def}{=} p^{\beta}, \eqno(6)
$$
but here $ \beta = \const > 0, \ p \in (1,\infty) $ so that
$$
\supp \psi(1,\infty; 0,-\beta;\cdot) = (1,\infty).
$$
 The space $ G\psi(1,\infty; 0, -\beta;\cdot)  $ coincides up to norm equivalence with
 the Orlicz space over the set  $ D $ with usually Lebesgue measure and with the
 correspondent $ N(\cdot) $ function

 $$
 N(u) = \exp\left(u^{1/\beta} \right), \ u \ge 1.
 $$
 Recall that the domain $ D $ has finite measure; therefore the behavior of the function
$ N(\cdot) $ is'nt essential.\par
{\bf Remark 1.} If we define the {\it degenerate } $ \psi_r(p), r = \const \ge 1 $ function as follows:
$$
\psi_r(p) = \infty, \ p \ne r; \psi_r(r) = 1
$$
and agree $ C/\infty = 0, C = \const > 0, $ then the $ G\psi_r(\cdot) $ space coincides
with the classical Lebesgue space $ L_r. $ \par
{\bf Remark 2.} Let $ \xi: D \to R $ be some (measurable) function from the set
$ L(p_1, p_2), \ 1 \le p_1 < p_2 \le \infty. $ We can introduce the so-called
{\it natural} choice $ \psi_{\xi}(p)$  as as follows:

$$
\psi_{\xi}(p) \stackrel{def}{=} |\xi|_p; \ p \in (p_1,p_2).
$$

\vspace{3mm}
{\bf B. Sobolev's-Grand Lebesgue Spaces.}\par
\vspace{3mm}
{\bf Definition.}\par
Let $ l=1,2,\ldots $ be any integer positive number. We introduce the following
so-called {\it Sobolev-Grand Lebesgue Space} $ W^l G(\psi), \ \psi \in \Psi(A,B) $
as a space of all weak $ l - $ times differentiable functions (in the Sobolev's sense)
with finite norm
$$
||f||W^l G(\psi) \stackrel{def}{=} ||f||G(\psi) + ||\nabla^l f||G(\psi). \eqno(7a)
$$
 Note that in the considered case, i.e. when $ f_{\partial D} = 0, $ the
$ W^l G(\psi) $ norm is equivalent the simple norm $ ||\nabla^l f||G(\psi): $

$$
||f||W^l G(\psi) \asymp ||\nabla^l f||G(\psi). \eqno(7b)
$$
 It is evident that the spaces $ W^l G(\psi), \ l \ge 1 $ are not rearrangement invariant.\par

\vspace{3mm}
{\bf C. Fundamental and truncated fundamental functions.}\par

\vspace{3mm}
Recall that if the rearrangement invariant space $ Y $ with the norm $ ||\cdot||Y $ over
the measurable space  $ (Z,\Sigma) $ equipped with the (non-trivial) measure $ \mu, $ then its fundamental function $ \phi(Y;\delta), \ \delta \in (0,\infty) $ is defined by
follows:
  $$
  \phi(Y;\delta) = \sup_{A, \mu(A) \le \delta} ||I(A)||Y.
  $$
 More detail information about the fundamental functions for rearrangement invariant
 spaces see in the book G.Bennet \cite{Bennet1}, chapter 3.\par
  The expression for the fundamental function for the Grand Lebesgue spaces  $ G(\psi) $
  may be written as follows:

  $$
  \phi(G(\psi),\delta) = \sup_{p \in (A,B)} \frac{\delta^{1/p}}{\psi(p)}. \eqno(8)
  $$
The fundamental function for the $ G\psi(A,B; \alpha,\beta; \cdot) $ spaces are calculated in the article \cite{Ostrovsky8}; see also \cite{Ostrovsky9}.\par
 We recall that when $ \beta > 0 $ as $ \delta \to 0+ $
$$
\phi(G\psi(a,b; \alpha,\beta); \delta)) \sim \beta^{-\beta} b^{-2\beta} \delta^{1/b} \
|\log \delta|^{-\beta}
$$
and

$$
\phi(G\psi(1,\infty; 0,-\beta); \delta)) \sim e^{-\beta} \ \beta^{\beta} \ |\log \delta|^{-\beta}, \ \delta \to 0+.
$$

{\bf Definition.}\par
We define the so-called {\it truncated fundamental function }
$ \phi_{p_-,p_+} (G(\psi); \delta) $ (only for GLS spaces)
as follows. Let $ p_- = \const \ge 1, \ p_+ = \const \in (p_-, \infty]. $ We put
$$
 \phi_{p_-,p_+}(G(\psi); \delta) \stackrel{def}{=}
\sup_{p \in (p_-,p_+) \cap \supp(\psi)} \frac{\delta^{1/p}}{\psi(p)},\eqno(9)
$$
where  also $ \delta = \const \in (0,\infty). $ \par
It is presumed that
$$
(p_-,p_+) \cap \supp(\psi) \ne \emptyset,
$$
as long as in the opposite case

$$
 \phi_{p_-,p_+}(G(\psi); \delta) = \infty
$$
and the formulating further main result, theorem 1, is trivial.\par
 It is evident that if $ \supp \psi \subset (p_-, p_+),  $ then
$$
\phi_{p_-,p_+}(G(\psi); \delta) = \phi(G(\psi); \delta).
$$

\vspace{3mm}

  \section{Main result: Sobolev's continuity inequality for Grand Lebesgue Spaces.}

\vspace{3mm}
We suppose in this section that $ d \ge 2 $ (the one-dimensional case $ d = 1 $ will be
investigated further) and that the (given) function $ f(\cdot), \ f_{\partial D}=0, $
belongs to some $ W^1 G\psi $ space

$$
f \in W^1 G\psi, \ \supp(\psi) = (A,B),
$$
where $ B > d; $ it may be considered also the case $ B = \infty. $ \par
 Let us denote
 $$
 A(1) = \max(A,d).
 $$

\vspace{3mm}
{\bf Theorem 1}. {\it The following Sobolev-type continuity inequality holds:}

$$
\omega(f, \delta) \le \frac{C(d,D) \ \delta \ ||\nabla f ||G(\psi)}
{\phi_{A(1),B}(G(\psi), \delta^d)}, \ \delta \in (0,1/e). \eqno(10)
$$

\vspace{3mm}

{\bf Remark 3.} In the when $ A \ge d $ the  last inequality may be
rewritten as follows:

$$
\omega(f, \delta) \le \frac{C(d,D) \ \delta \ ||\nabla f ||G(\psi)}
{\phi(G(\psi), \delta^d)}, \ \delta \in (0,1/e).
$$

{\bf Proof.} \par
Let $ f \in W^1 G\psi, \ f_{\partial D} = 0. $ We can and will assume without loss
of generality that

$$
||f||W^1 G(\psi) = 1,
$$
or equally

$$
|\nabla f|_p \le \psi(p), \ p \in (A,B).
$$

We have using the inequality (2c) for the values $ \delta \in (0,1/e) $ and
$ p \in (A(1),B): $

$$
\frac{\omega(f,\delta)}{\delta} \le C(d,D) \ \frac{p-1}{p-d} \ \delta^{-d/p} \ \psi(p)\le
$$

$$
\frac{C_2(d,D) }{\delta^{d/p}/ \{[(p-d)/(p-1)] \psi(p) \}}\le
\frac{C_3(d,D) }{\delta^{d/p}/\psi(p)}.
$$
 Since the last inequality is true for all the values $ p $ from the interval
 $ p \in (A(1),B), $ we obtain taking the minimum over $ p \in (A(1),B): $

$$
\frac{\omega(f,\delta)}{\delta} \le
\inf_{p \in (A(1),B)} \left[ \frac{C_3(d,D) }{\delta^{d/p}/\psi(p)} \right] =
$$

$$
 \frac{C_3(d,D) }{\sup_{p \in (A(1),B)} \left[\delta^{d/p}/\psi(p) \right]}=
\frac{C_3(d,D) }{\phi_{(A(1),B)}(G(\psi),\delta^{d})}.
$$
 This completes the proof of theorem 1. \par
\vspace{3mm}

\section{Examples.}

\vspace{3mm}
 We consider in this section some examples in order to show the exactness of
 the assertion of theorem 1.\par
 {\bf Theorem 2.} {\it For all the values } $ d = 2,3,4,\ldots $
 {\it there exists an admissible domain $ D,$ a function} $ \psi_0(\cdot) \in \Psi(1,\infty) $ {\it  and a non-trivial function } $ f_0(\cdot) \in G\psi_0,
 \ f_0:D \to R, $ {\it for which}

 $$
 \underline{\lim}_{\delta \to 0+}
 \left[ \omega(f_0,\delta): \frac{\delta}{\phi(G\psi_0,\delta^d)}  \right] > 0.\eqno(11)
 $$

 {\bf Proof.} Let us consider the space $ \phi(G\psi(1,\infty; 0,-\beta)) $ and a function
 $$
 f_0(x) = I(|x| \le 1) \ |x| \ |\log |x| \ |^{\beta}, \ \beta = \const > 0.
 $$
Note that the function $ f_0 $ is radial function, i.e. it dependent only on the Euclidean norm of a vector $ x, $ and that the function $ \psi(1,\infty;0,-\beta;p) $
asymptotically as $ p \to \infty $ coincides with the natural function for the function $ f_0. $ \par
 Here $ D = B. $ It is evident that as $  \delta \to 0+ $
$$
\omega(f_0, \delta) \sim \delta \ |\log \delta|^{\beta}.
$$
Recall that

$$
\phi(G\psi(1,\infty; 0,-\beta); \delta) \sim C_1(d,\beta) \ |\log \delta|^{-\beta}, \ \delta \to 0+.
$$
Further, we find by direct calculation as $ p \to \infty $ using the multidimensional
polar coordinates:

$$
|\nabla f_0|_p^p \sim \Omega(d)  \int_0^1 z^{d-1} \ |\log z|^{p \beta} \ dz =
$$

$$
= C_2(d) \int_0^{\infty} e^{-d y} \ y^{ dp} \ dy = C_2(d) d^{-dp+1} \Gamma(\beta p + 1);
$$

$$
|f_0|_p \sim C_3(d) \ \beta^{\beta} \ e^{-\beta} \ p^{\beta};
$$
we used Stirling's formula. \par
We conclude that

$$
f_0 \in G\psi_0, \ \psi_0(p) = \psi(1,\infty; 0,-\beta;p).
$$
 We obtain substituting into expression for $ \delta/\phi(G\psi_0,\delta)  $ as
 $ \delta  \to 0+: $

$$
\frac{\delta}{\phi(G\psi_0,\delta)} \sim C_4(d) \delta \ |\log \delta|^{\beta}.
$$
This completes the proof of theorem 1.\par
{\bf Remark 4.} Let us consider for comparison the case of the space
$ G\psi(A,B;\alpha,\beta; \cdot).$  Namely, we consider the following function

$$
g(x) = [\alpha/(\alpha - 1)] \ I(|x| \le 1) \ |x|^{-1/\alpha} \ |\log |x| \ |^{\gamma}.
$$
Here $ D = B \in R^d, \ d \ge 2, \ \alpha = \const > 1, \ \gamma = \const > 0; $ and denote

$$
b = \alpha d, \ \beta = \gamma + 1/b = \gamma + 1/(\alpha d);
$$
then $ d < b < \infty. $ \par
 We find by direct computation as $  p \to b-0: \ f \in L(1,b); $

$$
|\nabla g|_p^p \sim \Omega(d) \ \int_0^1 z^{d-1 - p/\alpha} \ |\log z|^{\gamma p} \ dz \sim
$$

$$
\Omega(d) \ \int_0^{\infty} e^{-y(d-p/\alpha)} y^{\gamma p} \ dy =
\Omega(d) \ \frac{\Gamma(\gamma p + 1)}{(d-p/\alpha)^{\gamma p + 1}}, \ p \in (1,b);
$$
therefore

$$
|\nabla g|_p \sim C_4(d,\alpha,\gamma) (b-p)^{-\gamma - 1/b} =
C_4(d,\alpha,\gamma) (b-p)^{-\beta}, \ p \in (1,b).
$$
On the other words, the function $  g(\cdot) $ belongs to the space
$  G\psi(1,b;0,\beta). $ \par
 It follows from the theorem 1 that

 $$
 \omega(g,\delta) \le C_5(d,\alpha,\gamma) \delta^{1-1/\alpha} \
 |\log \delta|^{\gamma + 1/b}, \ \delta \in (0, 1/e),
 $$
but really

$$
\omega(g,\delta) \sim C_6(d,\alpha,\gamma) \ \delta^{1-1/\alpha} \ |\log \delta|^{\gamma}.
$$
 Note that the main members in the two last expressions  coincides; but the second
members coincides only asymptotically, as $ \alpha d \to \infty. $ \par

\vspace{3mm}

\section{The one-dimensional case.}

\vspace{3mm}
 We consider in this section separately the one-dimensional case $ d = 1 $ and
 correspondingly the set $ D = [0,1], $ as long as we can obtain in the considered
 case the {\it asymptotical exact } as $ \delta \to 0+ $ value of an embedding constants.\par
 We suppose as before that $ f(0) = f(1) = 0 $ and that $ |\nabla f| \in G\psi, \
 \psi \in \Psi(1,\infty). $ \par

\vspace{3mm}
 {\bf Theorem 3.}

 $$
 \omega(f,\delta) \le 1 \cdot \frac{\delta \ ||f||G\psi}{\phi(G\psi,\delta) },
 \eqno(13)
 $$
 {\it when the constant "1" is the best possible.} \par

\vspace{3mm}
{\bf 1.} We obtain first of all the {\it upper bound} for Sobolev-Grand Lebesgue
continuity inequality in the one-dimensional case. Namely, let $ f(0) = f(1) = 0 $ and
$ \nabla f \in G\psi, \ \psi(\cdot) \in \Psi(A,B), $ i.e.

$$
|f^/|_p \le ||f^/||G(\psi) \cdot \psi(p), \ p \in (A,B), \ 1 \le A < B \le \infty.
$$
As long as
$$
f(y) - f(x) = \int_x^y f^/(z) dz, 0 \le x \le y \le 1,
$$
we have denoting $ \delta = |y-x| $ and using H\"older inequality:
$$
|f(y) - f(x)| \le |y-x|^{1-1/p} \ |f^/|_p \le ||f^/||G(\psi) \cdot \psi(p) \ \cdot
\delta^{1-1/p};
$$

$$
\omega(f,\delta) \le \delta \ ||f||G(\psi) \ \delta^{-1/p} \ \psi(p),
$$
therefore

$$
\omega(f,\delta) \le \delta \ ||f||G(\psi) \ \inf_{p \in (A,B)}
\left[\delta^{-1/p} \ \psi(p) \right] =
$$

$$
\omega(f,\delta) \le \delta \ ||f||G(\psi) \cdot
\frac{1}{ \sup_{p \in (A,B)}
\left[\delta^{1/p} / \psi(p) \right]} =
\frac{\delta ||f||G(\psi)}{\phi(G(\psi),\delta)}.
$$

\vspace{3mm}
{\bf 2.} Let us prove that the last inequality is in general case asymptotically
as $ \delta \to 0+ $ exact. Namely, we consider the following example (more exactly,
the family of examples) of a view:

$$
f_{\Delta}(x) = I( x \in [0,1]) \ x \ |\log x|^{\Delta}, \ \Delta = \const > 0.\eqno(14)
$$
 It is evident that as $ \delta \to 0+ $

 $$
 \omega(f_{\Delta}, \delta)  \sim \delta \ |\log \delta|^{\Delta};
 $$

$$
|\nabla f_{\Delta}|_p \sim \Delta^{\Delta} \ e^{-\Delta} \ p^{\Delta}, \ p \to \infty,
$$
and we choose as before
$$
\psi_{\Delta}(p) = |f_{\Delta}|_p;
$$
then

$$
||f_{\Delta}||G\psi_{\Delta} = 1.
$$
Further,

$$
\phi(G\psi_{\Delta}, \delta) \sim \sup_{p \in (1,\infty)} \frac{\delta^{1/p}}
{p^{\Delta} \Delta^{\Delta} e^{-\Delta} } \sim |\log \delta|^{-\Delta}.
$$
Thus,

$$
\underline{\lim}_{\delta \to 0+}
\left[\omega(f_{\Delta}, \delta):\frac{\delta ||f_{\Delta}||G\psi_{\Delta}}
{\phi(G\psi_{\Delta}, \delta)}  \right] =
\lim_{\delta \to 0+} \frac{\delta \ |\log \delta|^{\Delta}}
{\delta \ |\log \delta|^{\Delta}} = 1,
$$
Q.E.D.\par

\vspace{3mm}

\section{Hight derivatives.}

\vspace{3mm}
Let $ k,l $ be any positive integer numbers such that $ l - k \ge 1; $   and
$ \psi \in \Psi(A,B), $ where $ B > d/(l-k). $  We denote

$$
p(1) = d/(l-k), \ p(2) = d/(l-k-1); \ d/0 \stackrel{def}{=} + \infty;
$$

$$
(A(3), B(3)) = [(p(1), p(2))]\cap[(A(3), B(3))]
$$
and assume that $ (A(3), B(3)) \ne \emptyset. $\par
In this section we suppose for simplicity
$$
f_{\partial D}=0, \nabla f_{\partial D}=0, \ldots, \nabla^{l-1}f_{\partial D}=0.
$$

\vspace{3mm}

{\bf Theorem 4}. {\it The following generalized  Sobolev-Grand Lebesgue Space inequality
holds:}
$$
\omega \left(\nabla^k f, \delta \right) \le
\frac{C(d;l,k;D) \ \delta^{l-k} \ ||\nabla^l f||G\psi }
{\phi_{A(3), B(3)}(G(\psi), \delta^d)}.\eqno(15)
$$

\vspace{3mm}
{\bf Proof.} Let $ \nabla^l f \in G\psi, \ \psi(\cdot) \in \Psi(A,B), $ or equally

$$
| \nabla^l f |_p  \le || \nabla^l f ||G\psi \cdot \psi(p), \ p \in (A,B).
$$

We will use the following Sobolev's continuity inequality for the classical $ L_p $
spaces, see, e.g., \cite{Maz'ja1}, chapter 1, p. 60-64:

$$
\frac{|\nabla^k f(x) - \nabla^k f(y)|}{|x - y|^{\lambda}} \le C_6^{-1}(k,l;d,D;p) \
|\nabla^l f|_p.
$$
Here the constant $ C_6(\cdot) $ is bounded in the interval $ p \in (p(1), p(2)), $

$$
\lambda = l - k - d/p, \ (l-k-1)p < d < (l - k)p
$$
or equally $ p \in (p(1), p(2)). $ \par
 The last inequality may be rewritten (under our notations and conditions) as follows:

 $$
 \omega \left(\nabla^k f, \delta \right) \le C_6^{-1}(\cdot) \ \delta^{l-k-d/p}
  \cdot \psi(p) \cdot ||\nabla^l f||G\psi.\eqno(16)
 $$
The assertion of theorem 4 may be obtained as the proof of theorems 1 and 3 after the
dividing over $ \delta^{l-k} $ and taking minimum over $ p \in (A(3), B(3)). $ \par

\vspace{3mm}

\section{Concluding remarks. Generalizations.}

\vspace{3mm}
{\bf 1.} Let us denote

$$
\eta(\delta) = \frac{\delta \ ||\nabla f||G\psi }
{\phi_{A(1), B}(G(\psi), \delta^d)}, \ \delta \in (0,1/e),
$$
and introduce the generalized H\"older space $ H(\eta) $ as a space of continuous a.e. functions with zero boundary values $ f: D \to R $  with finite norm

$$
||f||H(\eta) \stackrel{def}{=} \sup_{x \in D} |f(x)| + \sup_{\delta \in (0,1/e)}
\left[\frac{\omega(f,\delta)}{\eta(\delta)} \right].
$$
 Then the assertion of theorem 1 may be reformulated as an continuous embedding
 theorem $ W^1 G\psi \subset H(\eta): $

$$
||f||H(\eta) \le C \ ||\nabla f||G\psi. \eqno(17)
$$

\vspace{3mm}
{\bf 2.} At the same examples as in the section 4 are true in the case when

$$
\psi(p) = \psi_L(p) \stackrel{def}{=} p^{\beta} \ L(p), \ p \in (1,\infty),
$$
or
$$
\psi^{(L)}(p) \stackrel{def}{=} (b-p)^{-\beta} \ L(1/(b-p)), \ p \in (1,b), b =
\const > d,
$$
where $ L=L(u) $ is continuous positive slowly varying as $ u \to \infty $ function.\par
The corresponding examples of the functions $ \{f = f(x) \} $ for the case when
$ D = B \subset R^d $ are described in \cite{Ostrovsky8}; see also \cite{Ostrovsky9}. \par
For instance, in the case when $ \psi(p) = \psi^{(L)}(p) $ the example function
$ f = f(x) $ has a view
$$
f(x) = |x| \ |\log |x| \ |^{\beta} \ L(1 + |\log |x| \ |) \ I(|x| \le 1).
$$

\vspace{3mm}
{\bf 3.} Some slight generalizations.\par
 Let now $ D, \ D \subset R^d $  be arbitrary open  domain in the space $ R^d. $
 We denote for arbitrary subset $ K $ of the region $ D, \ K \subset D $ by $ c_p(K) $ the $ p - $ {\it conductivity} of the set $ K; $ see the book of Maz'ja \cite{Maz'ja1}, chapter 4, section 4.1, p, 191-194 for the definition and some properties of this notion.  \par
  Introduce also as in \cite{Maz'ja1}, chapter 5, sections 5.3-5.4 the following functions:

 $$
 \gamma_p(x,y) = c_p[(D \setminus x) \setminus y]^{-1/p},
 $$

$$
\Lambda_p(\delta) = \sup_{x,y \in D, |x-y| \le \delta} \gamma_p(x,y),
$$

$$
\lambda^{(\psi)}(\delta) \stackrel{def}{=} \inf_{p \in (A,B)}
\left[ \Lambda_p(\delta) \ \psi(p) \right].
$$

\vspace{3mm}
{\bf Theorem 5.}\par
\vspace{3mm}
 $$
 \omega(f,\delta) \le \lambda^{(\psi)}(\delta) \cdot ||\nabla f||G\psi. \eqno(18)
 $$
{\bf Proof.} It is proved in  \cite{Maz'ja1}, chapter 5, sections 5.3-5.4 that
$$
\omega(f,\delta) \le \Lambda_p(\delta) \ |\nabla f|_p.
$$
 Therefore, if $ |\nabla f| \in G(\psi), \ \exists \psi \in \Psi(A,B), $ then

 $$
 \omega(f,\delta) \le \inf_{p \in \supp \psi} \left[\Lambda_p(\delta) \
  |\nabla f|_p \right]\le
 $$

 $$
\inf_{p \in \supp \psi} \left[\Lambda_p(\delta) \ \psi(p) \ ||\nabla f||G(\psi) \right] =
  \lambda^{(\psi)}(\delta) \cdot ||\nabla f||G\psi.
 $$
{\bf Remark 5.} The last result may be used, e.g., for the domains $ \{ D \} $ with
complicated boundaries.\par

\vspace{3mm}

{\bf 4. } Non-compactness of an embedding operator. \par
Let $ b = \const > 1, \ \beta = \const > 0, $

$$
\psi_{b,\beta}(p) = \psi(1,b;0,\beta + 1/b;p) = (b-p)^{-\beta - 1/b}, \ p \in (1,b),
$$

$$
\eta_{b,\beta}(\delta) = I(0 \le |x| \le 1) \ \delta^{1-1/b} \ |\log x|^{\beta}.
$$

 Let us denote also by $ E $ the unit embedding operator from the space
 $ W^1 G\psi_{b,\beta} $ into the space $ H(\eta_{b,\beta}): $

$$
E u = v, \ u \in W^1 G\psi_{b,\beta}, \ v \in H(\eta_{b,\beta}), \ u = v.
$$

{\bf Theorem 6.} {\it The operator} $ E $ {\it is'nt compact operator.} \par
{\bf Proof.} It is sufficient to consider only the one-dimensional case $ d = 1, $
i.e. $ D = [0,1]. $ \par
Let us consider the function

$$
g(x) =  I(0 \le x \le 1)  \ x^{-1/b} \ |\log x|^{\beta},
$$
an introduce the a {\it family} of a shift functions

$$
g_h(x) = T_h g(x) = g(x+h), \ x+h \le 1; \ T_h g(x) = g(x+h-1), x+h > 1.
$$
Here $ h \in (0, 1/2). $  It is evident that for  both the norms $ W^1 G\psi_{b,\beta} $
and $ H(\eta_{b,\beta} ) $

$$
||g_h||W^1 G\psi_{b,\beta} = ||g||W^1 G\psi_{b,\beta},
$$

$$
||g_h||H(\eta_{b,\beta}) = ||g||H(\eta_{b,\beta}),
$$
i.e. both the expressions does not dependent on the variable $ h. $\par
Let us calculate at first the norm $ ||g||W^1 G\psi_{b,\beta}. $  We have as
$ p \to b - 0: $
$$
|g_h|_p^p = |g|_p^p \sim \int_0^1 x^{-p/b} \ |\log x|^{\beta p} \ dx  =
b^{\beta p + 1} \frac{\Gamma(\beta p + 1)}{(b-p)^{\beta p + 1}};
$$

$$
|g_h|_p \sim b^{\beta + 1/b} \frac{\Gamma^{1/b}(\beta b + 1)}{(b-p)^{\beta + 1/b}}.
$$
Therefore, the family of the functions $ \{g_h\} $  belongs to some non-trivial ball
 in the space $ W^1 G\psi_{b,\beta}. $  Further,

 $$
 \omega(g_h,\delta) = \omega(g,\delta) = \eta_{b,\beta}(\delta), \ \delta \in (0,1/e).
 $$
This means that

$$
\sup_{h \in (0,1/2)} ||g_h||H(\eta_{b,\beta}) = 1.
$$
It is sufficient to prove that

$$
\overline{\lim}_{|h(1)-h(2)| \to 0} ||g_{h(1)} - g_{h(2)}||H(\eta_{b,\beta}) > 0,
$$
or equally

$$
\overline{\lim}_{h \to 0+}  \zeta(h) > 0,
$$

$$
\zeta(h) \stackrel{def}{=} ||g_{h} - g||H(\eta_{b,\beta}).
$$
We get:

$$
\zeta(h) \ge \sup_{\delta \in (0,1/e)} \frac{\omega(g_h-g,\delta)}
{\eta_{b,\beta}(\delta)} =
$$

$$
\sup_{\delta \in (0,1/e)} \sup_{|\tau| \le \delta} \sup_{x \in [0,1]}
\frac{|g(x+\tau + h) - g(x+h) - g(x+\tau) + g(x)|}{\eta_{b,\beta}(\delta)} \ge
$$

$$
\sup_{\delta \in (0,1/e)} \sup_{|\tau| \le \delta}
\frac{|g(\tau + h) - g(h) - g(\tau)|}{\eta_{b,\beta}(\delta)}.
$$

We conclude taking the values $ \tau = \delta = h $ that for all sufficiently small positive  values $ h $

$$
\zeta(h) \ge 2 - 2^{1-1/b} = \const > 0.
$$
 This completes the proof of theorem 6.

\vspace{3mm}

{\bf 5.} Note that the inequality (10) contains as a particular case the classical
result (2c) for ordinary Lebesgue spaces $ L_p, \ p \ge d, $ as long as the fundamental function for these  spaces has a view

$$
\phi(L_p, \delta) = \delta^{1/p}.
$$

\end{document}